\documentclass[a4paper,12pt]{amsart} 
\usepackage{amssymb,euscript,eufrak,verbatim}
\usepackage[utf8]{inputenc}      
\usepackage[T1]{fontenc}
\usepackage[polish,english]{babel}

\usepackage{tikz}
\usetikzlibrary{snakes}
\usetikzlibrary{arrows}

\def\a{\alpha} 

\def\z{\phi}
\def\e{\varepsilon}

\def\S1{\mathbb{S}^1}

\def\R{\mathbb{R}}

\def\B{\mathcal{B}}

\def\d{\delta}

\def\a{\vec a}

\newcommand{\vt}{{\vec\theta}}

\newcommand{\inner}[2]{\langle{#1},{#2}\rangle_\theta}
\newcommand{\hf}{{h_{\va}}}

\def\i{{i,j}}
\def\va{\vec{a}}
\def\vb{\vec{b}}

\newtheorem{theorem}{Theorem}

\DeclareMathSymbol{\varnothing}{\mathord}{AMSb}{"3F} 
\begin{document}
\renewenvironment{proof}{\noindent {\bf Proof.}}{ \hfill\qed\\ }
\newenvironment{proofof}[1]{\noindent {\bf Proof of #1.}}{ \hfill\qed\\ }

\title{Weak-mixing polygonal billiards}  

\author {Alba M\'alaga Sabogal}
\address{D\'epartement de Math\'ematiques,B\^atiment 425, Facult\'e des Sciences D'Orsay, Universit\'e Paris-Sud 11, F-91405 Orsay Cedex.}
\email{alba.malaga@polytechnique.edu}

\def\curraddrname{{\itshape Address}}
\author{Serge Troubetzkoy}
\address{Aix Marseille Universit\'e, CNRS, Centrale Marseille, I2M, UMR
  7373, 13453 Marseille, France}
  \curraddr{I2M, Luminy\\ Case 907\\ F-13288 Marseille CEDEX 9\\ France}

 \email{serge.troubetzkoy@univ-amu.fr}
\begin{abstract}  
We consider the set of polygons all of whose sides are vertical or horizontal with 
fixed combinatorics (for example all the figure "L"s). We show that there is 
a dense $G_\delta$ subset  of such polygons such that for each polygon in this $G_\delta$ set 
the billiard flow is weakly-mixing in almost every direction.
\end{abstract} 
\maketitle

\section{Introduction}
A polygon is called rational if all the angles between pairs of sides are rational multiples of $\pi$.  The billiard flow
in a rational polygon is quasi-intergrable. Namely if we fix the initial direction of the flow, then the billiard orbit is contained in a finite set of directions, this set of directions is the integral.
 it has a natural first integral, the direction of the flow.  A major achievement in the theory of
polygonal billiards is the Kerchoff, Masur, Smillie theorem: the billiard flow  in any rational polygon is ergodic, in fact uniquely ergodic, for almost
every value of the first integral \cite{KMS}.
One of the major unsolved questions about rational polygonal billiards is to understand if for almost all values of
the integral the billiard flow is in fact weakly mixing.

There is a well known construction, known as unfolding, which associates a translation surface to each
rational polygon.   \'Avila and Forni have  shown  that the typical translation flow on a typical surface is weakly mixing
\cite{AF}, but this
result tells us nothing about the behavior of polygonal billiards since the translation surface  formed by unfolding are a set of measure zero in the space of all translation surfaces.  Recently \'Avila and Delecroix have shown that the same holds for the billiard
in a regular polygon (except for the lattice cases: the equilateral triangle, the square, or the regular hexagon) \cite{AD}.  Other tham this, to our knowledge,
there is only one result about weak-mixing for the billiard flow in a rational polygon. Consider polygons for which
all sides are vertical or horizontal (VH-tables).  Fix the combinatorics of such a polygon (see the next section for an exact definition) and fix an
angle $\theta \in \mathbb{S}^1$, then Gutkin and Katok have shown that there  is 
 a dense $G_\delta$ set of VH-pollygons with this combinatorics, such that for each polygon in this dense $G_\delta$ set
 the billiard flow in the direction $\theta$ is weakly mining \cite{GK}. Their result also holds for tables of fixed combinatorics
 for which all sides are parallel to the sides of a fixed equilateral triangle.
 
 In this article we strengthen this last result.  We show that for any fixed combinatorics of VH-tables, there  is 
 a dense $G_\delta$ set of VH-pollygons with this combinatorics, such that for each polygon in this dense $G_\delta$ set
 the billiard flow is weakly mining is weakly-mixing in almost every direction. The method of proof is the same as in \cite{GK}, their
 article contains several of the necessary ingredients for our proof.

\section{Definitions and main results}
For details on polygonal billiards see \cite{MT}.
A \emph{VH-polygon} is a simply connected polygon all of whose sides are either vertical
or horizontal.  The smallest rectangle congaing  a VH-polygon is called is called its {\em bounding box}.  
We orient the border of $\B$ in such a way that the interior of the table is always
to the left of the border, thus
we  associate a finite  word in four letters $W$, $S$, $E$, $N$ to the VH-polygon $P$. This word is assumed to start at the side of the table on the
lower side of the bounding box, if there are several we take the left most segment. This word is called the \emph{combinatorics} of $P$. 

\begin{figure}[h]
\centering
\begin{tikzpicture}[scale=0.75]

\draw[thin] (2,1) -- node[below]{\footnotesize{$\ell_1$}} (3,1) -- (3,1.5) -- (3.5,1.5) -- (3.5,1) -- (4,1) -- node[right]{\footnotesize{$\ell_6$}}(4,2) -- (5,2) -- node[right]{\footnotesize{$\ell_8$}}(5,3) -- (3,3) -- (3,4) -- (2,4) -- (2,3) -- (1,3) --(1,2) --(2,2) --cycle;

\draw[dotted](1,1) rectangle +(4,3);

\end{tikzpicture}
\caption{A VH-polygon and its bounding rectangle.}\label{fig0}
\end{figure}
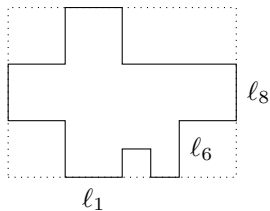

The location $(x,y)$ of the lower left corner of the bounding box, the combinatorics, and the lengths $(\ell_1,\dots,\ell_n)$ of sides of $P$ 
naturally parametrize the set of all VH-polygons.
The set of VH-polygons with a fixed combinatorics is naturally parametrized by an open subset  of $\R^{n}$. 
To see this let $\B$ be a table with a given combinatorics. Let $\ell_1,\dots,\ell_n$ be the lengths of the sides of $\B$. The total displacement to the right must be equal to the total displacement to the left, and the same for up and down. So there are two dependent length parameters for each word, but we need to also specify
the coordinates $(x,y)$ which compensate for this.

A \emph{VH-table} is a VH-polygon from whose interior we remove a finite number of disjoint VH-polygons.  
The \emph{combinatorics} of a VH-table is the collection of these words, with the word associated to the exterior component distinguished. We want to parametrized the set
of VH-tables with fixed combinatorics and fixed area for which the outer most bounding
box has lower left corner at the point $(1,1)$.
As before we consider the  lengths $(\ell_1,\dots,\ell_n)$
of all the sides of the table, the the set of VH-tables with a fixed combinatorics and 
fixed area is then naturally parametrized by an open subset $\mathcal O$ of $\mathbb{R}^{n-3}$ since
we additionally fixed the area and the point $(1,1)$. Some coordinates of
points in this space are lengths, others are positions of corners of bounding boxes, thus
we will use the letters $\vec a$ or $\vec b$ to denote points in $\mathcal O$. We will often confound the parameter
$\vec a$ with the table $\B^{\a}$ which it determines.

If all the entries of $\vec{a}$ are rational, the table can be tiled by rectangles of size
$(\frac1p,\frac1q)$ with $p,q$ positive integers. The number of rectangles $N(p,q)$ in
the tiling is $pq$ times the area of the table.
We call a table  in $\mathcal O$ a \emph{$(p,q)$-VH table} if it is tiled with rectangles of
size $(\frac1p,\frac1q)$ with
$p,q$ minimal.

Fix a VH-table $\a \in \mathcal{O}$ and  a direction $\theta \in \S1$.
The \emph{billiard flow $\phi_t^{\a,\theta}$ on $\B^{\a}$ in the direction $\theta$} or simply $\phi_t^\theta$  is  the free motion on the interior of $\B^{\a}$ with elastic collision from the boundary of $\B^{\a}$.
Once launched in the direction $\theta$, the billiard direction can only achieve four directions $[\theta] := \{\pm \theta, \pm(\pi - \theta)\}$;
thus the \emph{phase space} $X^{\a,\theta}$  
of the billiard flow in the direction $\theta$ is a subset of the cartesian product of $\B^{\a}$ with these four directions. 
Note that in this notation  $\phi_t^\theta$, $\phi_t^{-\theta}$, $\phi_t^{\pi-\theta}$ and $\phi_t^{\theta - \pi}$
are all the same.

For each direction $\theta$, the billiard flow  $\z^{\theta}_t$ preserves 
 the area measure $\mu$ on $\B^{\a}$ times a discrete measure on $[\theta]$, we
will also call this measure $\mu$. 
The billiard flow on the full phase space preserves the volume measure $\mu \times \lambda$ with $\lambda$ the length measure on $\S1$. 

A flow $\psi_t : M \to M$ preserving a measure $m$ is called 
{\em weakly-mixing} if for any $f \in L^2(M,m)$ the function
$t \to \langle U^tf,f \rangle$ strongly converges in the sense of Cesaro to $|\langle f,1 \rangle |^2$ as $|t| \to \infty$, 
where $1$ is the constant function taking the value 1 on $M$ and $U^tf= f \circ \psi_t$ \cite{P}.

 \begin{theorem}\label{GuKa}
 The set of VH-tables with a fixed area and fixed combinatorics (with at least 6 sides) contains a dense $G_\delta$ set of tables $G$ such that for each table $\vec a  \in G$ 
 there exists a full measure $G_\delta$-dense set $\Theta_{\vec a}$ of directions
 such that the billiard flow on the table $\vec a \in G$ is weakly mixing in every direction $\theta \in \Theta_{\vec a}$. 
 
In fact, the conclusion holds if we restrict ourselves to any closed subset of the set of VH-tables with fixed combinatorics which for each $Q \ge 1$ contains a dense set of $(p,q)$-VH tables with $\min(p,q) \ge Q$. \end{theorem}


%
 
\section{Proof of generic weak mixing}
\begin{proofof}{Theorem \ref{GuKa}} 
Let us consider $\va \in\mathcal O$. 
The  table $\B^{\va}$ can be unfolded \cite{MT} to a translation surface which in fact will be obtained as a gluing of four copies of $\B^{\va}$ along corresponding pairs of parallel sides.  Without loss of generality we assume that the area of the table is $\frac14$, and thus the
area of the unfolded surface is $1$.
For our proof we look at these copies as embedded disjointly in $\R^2$ in such a way that this embedding varies continuously with $f$. The bounding box of one copy of the table is always
placed with the inferior left corner to be $(1,1)$ and the other copies to be the images under reflections with respect to the axis, thus their respective bounding boxes having
a corner placed at $(\pm1,\pm1)$. 

Thus, modulo the boundary of $\B^{\va}$ which is a set of zero-measure, the phase space of the billiard on the table $\B^{\va}$ is a cartesian product of a compact subset of $\mathbb R^2$ with $\S1$. Thus any  function in $\mathcal L^2(\mathbb R^2\times\S1)$ is projected to a function in $\mathcal L^2(\B^{\va}\times\S1)$ by restriction and reciprocally, any function in $\mathcal L^2(\B^{\va}\times\S1)$ can be seen as a function in $\mathcal L^2(\mathbb R^2\times\S1)$. 

We extend the flow $\phi^{\va}_{t}$ to all of $\R^2$ by setting $\phi^{\va}_{t}(\vec z,\vt)=(\vec z,\vt)$ on the exterior of $\B^{\va}$.
For each $t$ let $U^{\va}_t$ be the unitary operator on $\mathcal L^2(\mathbb R^2\times\S1)$ defined by $U^{\va}_th(\vec z,\vt)=h(\phi^{\va,\theta}_t(\vec z,\vt))$.

Suppose additionally that $\va$ is a $(p,q)$-VH table.
Let $\chi^{\va}$ be the characteristic function of $\B^{\va}\times \S1$, and for any $h:\R^2\times\S1\to \R$ let $h_{\va}:=\chi^{\va}\cdot h$ be the restriction of $h$ to $\B^{\va}\times \S1$. 
Let $h^d : \R^2\times\S1 \to \R$ be the function defined by
$$\frac{1}{N(p,q)} \sum_{\i} \hf(s +(\i),\theta).$$
It follows from definition that 
\begin{equation}\label{e0}
\int_{\B^{\va}\times \S1}h^d\,d(\mu\times\lambda)=\int_{\B^{\va}\times \S1}h\,d(\mu\times\lambda).
\end{equation} Note that $\big(\hf)_{\va}=\hf$, and thus $\hf^d=h^d$. Furthermore ${h^d}_{\va}\ne \hf^d$ however ${h^d}_{\va}=\big(\hf^d\big)_{\va}$.
Note that $h^d$ is $(\frac1p,\frac1q)$-periodic on $\R^2\times  \{\theta\}$, for any $\theta$. 
Thus we can think of $h^d$ as descending to a function in $L^2([0,\frac1p] \times [0,\frac1q] \times \S1)$,
we will use the same symbol $h^d$ for both functions, but will differentiate
by noting the variable $(\vec z,\vt) \in \R^2 \times \S1$, and $(u,\theta) \in [0,\frac1p] \times [0,\frac1q] \times \S1$.
 By Theorem 3 of \cite{GK}, if $\theta$ is irrational, i.e. non commensurable with $\pi$, then
the function ${h^d}_{\va}$ is the projection of $\hf$ on the discrete spectrum of the
flow $\phi_t$. Let $$h^c:= h - h^d.$$ Then 
$$\hf={h^d}_{\va}+{h^c}_{\va},$$ thus ${h^c}_{\va}$ is the projection of $\hf$ on the continuous
spectrum of $\phi_t$.


For any $h_1,h_2$ in $\mathcal L^2(\R^2\times \S1)$, let  
$$\inner{h_1}{h_2}=\int_{\mathbb R^2}h_1(\vec z,\vt)h_2(\vec z,\vt)\,d\mu$$
and
$$|| h _1|| =\sqrt{ \inner{h_1}{h_1}}.$$
For any $h\in\mathcal L^2(\R^2)$, and any fixed direction $\theta$,
\begin{equation}\label{e1}
\begin{array}{l}
\inner{U^{\va}_t\hf}{\hf}-\big|\inner{\hf}{\chi^{\va}}\big|^2=\\
\inner{U^{\va}_t\big[({h^d}_{\va}-\inner{\hf}{\chi^{\va}}\chi^{\va})+\inner{\hf}{\chi^{\va}}\chi^{\va}+{h^c}_{\va}\big]}{\hf}-\big|\inner{\hf}{\chi^{\va}}\big|^2=\\
\inner{U^{\va}_t\big[{h^d}_{\va}-\inner{\hf}{\chi^{\va}}\chi^{\va}\big]}{\hf}+\inner{U^{\va}_t{h^c}_{\va}}{\hf}=\\
\inner{U^{\va}_t\big[{h^d}_{\va}-\inner{\hf}{\chi^{\va}}\chi^{\va}\big]}{\hf}+\inner{U^{\va}_t{h^c}_{\va}}{{h^c}_{\va}}.
\end{array}
\end{equation}

By the definition of the continuous spectrum for any $\e > 0$ there is a set $T \subset \R$ of density
one such that if $t \in T$ then 
\begin{equation}\label{e2}
\inner{U^{\va}_t{h^c}_{\va}}{{h^c}_{\va}} < \e.
\end{equation}

Suppose that $h$ is additionally continuous;  thus for any $\e >0$ there is a 
$\d(\e,h) > 0$ such that
$|h(\vec z,\vt) - h(s',\theta)| < \e/||h||$ if $|s-s'| < \d(\e,h)$.
Suppose that $\max(\frac1p,\frac1q) < \d(\e,h)$. Then
\begin{equation}\label{e3}
\begin{array}{rl}
|h^d(\vec z,\vt) - h^d(s',\theta)| =& |h^d(u,\theta) - h^d(u',\theta)|\\
\le & \frac1{pq} \sum_{\i} |h(u+(\i),\theta) - h(u'+(\i),\theta)|\\
< &\frac{\e}{||h||pq} \times \text{the number of rectangular tiles}\\
 \le & \frac{\e}{||h||}.
\end{array}
\end{equation}
 
Using the Cauchy-Schwartz inequality, and then applying Formulas \eqref{e0} and \eqref{e3} for any $t \in \R$ we have
\begin{equation}\label{e4}
\begin{array}{rl}
\left |\inner{U^{\va}_t\big[{h^d}_{\va}-\inner{\hf}{\chi^{\va}}\chi^{\va}\big]}{\hf}  \right |  \le &
\left | \left | U^{\va}_t [{h^d}_{\va}-\inner{\hf}{\chi^{\va}}\chi^{\va}] \right | \right |  \cdot
 \left |  \left | \hf \right | \right |\\
= &||  {h^d}_{\va}-\inner{\hf}{\chi^{\va}}\chi^{\va} ||  \cdot \left | \left |\hf  \right | \right |  \\
= &||  {h^d}_{\va}-\inner{{h^d}_{\va}}{\chi^{\va}}\chi^{\va} ||  \cdot \left | \left |\hf  \right | \right |  \\
< & \frac{\e}{||h||} \cdot  \left | \left |\hf  \right | \right | \\
\le & \e.
\end{array}
\end{equation}

Combining Equations \eqref{e1},\eqref{e2} and \eqref{e4} we conclude that for any $t \in T$
$$
\inner{U^{\va}_th_{\va}}{h_{\va}}-\big|\inner{h_{\va}}{\chi^{\va}}\big|^2 < 2\e.
$$

For a given VH-table $\vb$ in $\mathcal O$ and every $\tau>N$, we consider the set of directions

$$\Theta(\vb,\tau,h,N):=\left\{\theta:\exists t\in(N,\tau),\big|\inner{U^{\vb}_th_{\vb}}{h_{\vb}}-|\inner{h_{\vb}}{\chi^{\vb}}|^2\big|< \textstyle\frac1N\right\}.$$

From Formula \eqref{e2} it follows that for any  any $(p,q)$-VH-table $\va$ and any $\theta$ irrational there exists $\tau_0>0$ such that $\theta$ is in $\Theta(\va,\tau_0,h,N)$. Moreover, $\theta$ is in $\Theta(\va,\tau,h,N)$ for any $\tau\ge\tau_0$ as these sets are increasing in $\tau$. Thus, the measure of $\Theta(\va,\tau,h,N)$ tends to 1 as $\tau$ goes to infinity.

For any $\theta\in \Theta(\va,\tau,h,N)$, let $\eta_\theta>0$ be as large as possible such that the $\eta_\theta$-neighborhood $\mathcal U(\va,\eta_\theta)$ is included in 
$$\left\{ \vb \in\mathcal O:\exists t<\tau,\,\big|\inner{U^{\vb}_th_{\vb}}{h_{\vb}}-|\inner{h_{\vb}}{\chi^{\vb}}|^2\big|\le \textstyle\frac1{2N}\right\}$$

So, for any $N$ we can choose $\tau(\va,h,N)$ and $\eta(\va,h,N)>0$ such that $\overline\Theta(\va,h,N):=\Theta(\va,\tau(\va,h,N),h,N)\cap\{\theta:\eta_\theta>\eta(\va,h,N)\}$ has measure at least $1- \textstyle\frac1{N^2}$.

Thus, any $\vb$ in  $\eta$-neighborhood $\mathcal U(\va,\eta(\va,h,N))$ satisfies
$$\big|\inner{U^{\vb}_th_{\vb}}{h_{\vb}}-|\inner{h_{\vb}}{\chi^{\vb}}|^2\big|< \frac1N$$
for any $\theta$ in $\overline\Theta(\va,h,N)$ and some $N<t<\tau(\va,h,N)$. 

Let $\{h_j\}$ be a countable basis of continuous functions in $\mathcal L^2(\R^2)$. Under the identification $h_j(\vec z,\vt)=h_j(s)$, this is a countable basis of $\mathcal L^2(\R^2\times\{\theta\})$ for all $\theta$.
Let $\{\va_i\}$ be a dense collection of $(p,q)$-VH-tables. Let us consider the following dense $G_\delta$-set
$$
\mathcal G := \bigcap_{N=1}^\infty\bigcup_{i=N}^\infty\bigcap_{j=1}^N\mathcal U\big(\va_i,\eta(\va_i,h_j,N)\big).
$$

Let $\vb \in\mathcal G$. Then for each $N$ there is $i_N\ge N$ such that $\vb$ and $\va_{i_N}$ are $\eta_N$-close, where $\eta_N:=\min\big\{\eta(\va_{i_N},h_j,N):1\le j\le N\big\}$.

The set $\bigcap_{j=1}^N\overline\Theta(\va_{i_N},h_j,N)$ is of measure at least $1-\frac1N$ and for every direction $\theta$ in this set, and every $j\le N$, there exists $t_N(j)\ge N$ such that $$\big|\inner{U^{\vb}_{t_N(j)}h_{\vb}}{h_{\vb}}-|\inner{h_{\vb}}{\chi^{\vb}}|^2\big|< \frac1{N},$$ where $h=h_j$.

Let  $\displaystyle\hat\Theta_{\vb}:=\bigcap_{n=1}^\infty\bigcup_{N=n}^\infty \bigcap_{j=1}^N\overline\Theta(\va_{i_N},h_j,N)$. This is a dense $G_\delta$-set of full measure. For every $\theta\in\hat\Theta_{\vb}$ and any $j\ge 1$ there exist  sequences  $N_k=N_k(j,\theta,g)\to \infty$ and $ t_k=t_k(j,\theta,g)\ge N_k$,  such that 
\begin{equation}\label{e5}
\big|\inner{U^{\vb}_{t_k}h_{\vb}}{h_{\vb}}-|\inner{h_{\vb}}{\chi^{\vb}}|^2\big|< \frac1{N_k}
\end{equation}
where $h=h_j$. 

We will now show that for every $\theta \in\hat\Theta_{\vb}$ the flow $\phi_t^\theta$ on $\mathcal \B^{\vb}$ is weak-mixing. Fix $\theta\in\hat\Theta_{\vb}$. Since $\theta$ is fixed, we will denote by $U^{\vb}_t$ the operator $U^{\vb}_t$ restricted to the direction $\theta$.

Suppose $\phi^\theta_t$ is not weak-mixing, then there exists an eigenfunction $h\in\mathcal L^2(\B^{\vb}\times \{\theta\})$ with an eigenvalue $e^{i a}$; thus, $U^{\vb}_th=e^{i a t}h$. We can assume that : $||h||_\theta=1$  and $\int_{\B^{\vb}}h\,d\mu=0$.  

Let $\varepsilon>0$, and let $h_j$ be such that $||h-h_j||_\theta<\varepsilon$.  Then for any $t\in\R$ we have 
$$\begin{array}{rcl}U^{\vb}_th_j&=&U^{\vb}_t(h_j-h)+U^{\vb}_th\\ &=& U^{\vb}_t(h_j-h)+e^{i a t}h.\end{array}$$
Therefore
$$\begin{array}{rcl}\inner{U^{\vb}_th_j}{({h_j})_{\vb}}&=&\inner{U^{\vb}_t(h_j-h)+U^{\vb}_th}{({h_j})_{\vb}}\\
&=&\inner{U^{\vb}_t(h_j-h)+e^{i a t}h}{\big(({h_j})_{\vb}-h\big)+h}\\
&=&\inner{U^{\vb}_t(h_j-h)}{({h_j})_{\vb}-h}+2\Re\inner{h}{({h_j})_{\vb}-h}+e^{i a t}\end{array}$$
which implies the estimate
$$\left|\inner{U^{\vb}_th_j}{({h_j})_{\vb}}-e^{i a t}\right|<\varepsilon^2+2\varepsilon.$$

Since 
$$\begin{array}{rcl}\big|\int_{\B^{\vb}}{({h_j})_{\vb}}\,d\mu\big|&=&\big|\int_{\B^{\vb}}{(({h_j})_{\vb}-h)\,d\mu\big|} \\
&\le&||({h_j})_{\vb}-h||_\theta ||\chi^{\vb}||_\theta \;\text{ (by Cauchy-Schwartz inequality)}\\
&=&||({h_j})_{\vb}-h||_\theta \,\mu(\B^{\vb})=||({h_j})_{\vb}-h||_\theta\\
&<&\varepsilon,\end{array}$$ 
we have for any $t$ (by the triangular inequality):

$$\begin{array}{rcl}
\left|\inner{U^{\vb}_th_j}{({h_j})_{\vb}}-|\inner{h_j}{\chi^{\vb}}|^2\right|&\ge&|\inner{U^{\vb}_th_j}{({h_j})_{\vb}}|-|\inner{h_j}{\chi^{\vb}}|^2\\ 
&>&(|e^{i a t}|-\varepsilon^2-2\varepsilon)-\varepsilon^2\\
&>&1-2\varepsilon-2\varepsilon^2.
\end{array}$$

Thus, for $\varepsilon$ small enough, $h_j$ will satisfy
$$\left|\inner{U^{\vb}_th_j}{({h_j})_{\vb}}-|\inner{h_j}{\chi^{\vb}}|^2\right|>\textstyle\frac12$$
for all $t$. This contradicts Equation \eqref{e5}. Thus such an h does not exist and the direction $\theta$ is weak-mixing.
\end{proofof}

\section{Generalizations}
\begin{figure}[h]
\centering
\begin{minipage}[ht]{0.45\linewidth}
\centering
\begin{tikzpicture}[scale=0.75]
\draw (0,0) -- (1,1.732) -- (4,1.732) -- (5,0) -- (0,0);
\end{tikzpicture}
\end{minipage}
\begin{minipage}[ht]{0.45\linewidth}
\centering
\begin{tikzpicture}[scale=0.75]
\draw (0,0) -- (1,1.732) -- (5,1.732) -- (4,0) -- (0,0);
\end{tikzpicture}
\end{minipage}
\caption{}\label{fig2}
\end{figure}
The theorem also holds for billiards in polygons with  fixed combinatorics (and at least 4 sides) all of
whose sides are parallel to a fixed equilateral triangle.  In Figure \ref{fig2} we give two examples of fixed combinatoric:
there is a dense $G_\delta$ of  trapezoids with angles $\pi/3$ and $2\pi/3$ (of fixed area) for which the billiard flow is
weakly mixing in almost every direction, and the same for parallelograms of fixed area with these angles.

Our method also proves topological genericity of the set of translation surfaces within a given stratum for which the
translation flow is weakly mixing in almost every direction; the result in \cite{AF} only showed the measure theoretic
genericity of such translation surfaces. Indeed, it is sufficient to check that in each stratum the set of rectangle tiled surfaces
which are $(p,q)$-tiled with $\min(p,q) \ge Q$ is dense for each $Q$. 

There are two well known coordinate systems on the set of of translation surfaces: local coordinates for the ablian differential or
the zippered rectangle representation; the density holds for both.

Consider first local coordinates for Abelian differentials, they are given by the relative periods of the holomorpic one form, i.e.\
a collection of planar vectors corresponding to certain saddle connections (see \cite{Z} section 6.4).
Lemma 18 and Remark 7 of \cite{KZ} show that the set of compact surfaces
with single-cylinder directions is dense in each stratum.
The constructed surfaces have the vertical direction as single-cylinder.
By making arbitrarily small changes to the vectors of the vertical saddle connections,
we ensure that the lengths of all vertical saddle connections are rational
multiples of the width of the cylinder without changing the fact that the vertical
flow is a single-cylinder direction. By rescaling the width by an arbitrarily small
amount, we may ensure all such lengths are in fact rational. We then need only
translate all the charts of M by an arbitrarily small amount to ensure that all
singularities project to rational coordinates; the surface is now square-tiled with a
vertical single-cylinder direction.  If additionally we require that all these rational number have a fixed prime denominator
$q$, then the surface is $(q,q)$ tiled by squares, and is not square tiled by larger squares.  For any $Q \ge 1$ the set of 
surfaces of this type with $q \ge Q$ is dense.

Next we  describe the density in the zippered rectangle
representation of translation surfaces (see \cite{Z} sections 5.5 and 5.7 for a description of this representation). 
The zippered rectangle coordinates are given by an interval exchange,
thus a permutation and some lengths of intervals, the heights of the rectangles, as well as the heights of the singularities.
If each of the length and height parameters is
rational with common prime denominator $q$ then the surface is $(q,q)$- tiled by square  and is
not square tiled by larger squares.  For any $Q \ge 1$ the set of surfaces of this type with $q \ge Q$ is dense.


\begin{thebibliography}{99} \footnotesize{

\bibitem[AD]{AD} A.\ \'Avila and V.\ Delecroix \textit{Weak mixing directions in non-arithmetic Veech surfaces}
ArXiv 1304.3318v2.

\bibitem[AF]{AF} A.\ \'Avila and G.\ Forni \textit{ Weak mixing for interval exchange transformations and translation flows}
Annals of Math.\ (2) 165  (2007), no.\ 2 637-664. 


%


\bibitem[GK]{GK} E.\ Gutkin and A.\ Katok,
{\em Weakly mixing billiards} in Holomorphic dynamics (Mexico, 1986) 163--166,
Lecture Notes in Math., 1345, Springer Berlin, 1988.

\bibitem[KMS]{KMS} S.\ Kerckhoff, H.\ Masur and J.\ Smillie, {\em Ergodicity of billiard flows and quadratic differentials},  Annals of Math.\
(2) 124 (1986), no.\ 2, 293--311.

\bibitem[KZ]{KZ} M.\ Kontsevich and A.\ Zorich, {\em Connected components of the moduli spaces of abelian
differentials with prescribed singularities}, Inventiones mathematicae 153(3) (2003) 631--678.


\bibitem[MT]{MT} H.\ Masur and S.\ Tabachnikov, {\em  Rational billiards and
flat structures},  Handbook of dynamical systems, Vol.~1A, 1015--1089,
North-Holland, Amsterdam, 2002. 


\bibitem[P]{P} K.\ Petersen, {\em Ergodic theory} Cambridge University Press (1983).

\bibitem[Z]{Z} A.\ Zorich, \textit{Flat surfaces}, in Frontiers in number theory, physics, and geometry. {I}
437--583,  Springer, 2006.


}

\end{thebibliography}
\end{document}